%% file: main_JSV.tex
\documentclass{article}
\usepackage[utf8]{inputenc}
\usepackage[margin=1in]{geometry}

\usepackage{enumerate}
\usepackage{amsmath}
\usepackage{amssymb}
\usepackage{gensymb}
\usepackage[numbers,sort&compress]{natbib}
\usepackage{xcolor}
\usepackage{graphicx}
\usepackage{courier}

\usepackage{blindtext}
\newcommand\blfootnote[1]{%
\begingroup
\renewcommand\thefootnote{}\footnote{#1}%
\addtocounter{footnote}{-1}%
\endgroup
}

\newcommand{\bluenote}[1]{{#1}}
\newcommand{\bbmat}{\begin{bmatrix}}
\newcommand{\ebmat}{\end{bmatrix}}
\newcommand{\tp}{{\text{T}}}
\newcommand{\pCite}[1]{(\cite{#1})}
\newcommand{\FI}{{\text{FI}}}

\begin{document}

\title{Friction Induced Instability and Vibration in a Precision Motion Stage with a Friction Isolator}                      

\author{Jiamin Wang{$^1$}, Xin Dong{$^2$}, Oumar R. Barry{$^{1, 3}$}, and Chinedum Okwudire{$^2$}}
\date{}




\maketitle

\begin{abstract}
Motion stages are widely used for precision positioning in manufacturing and metrology applications. However, they suffer from nonlinear pre-motion (i.e., "static") friction, which adversely affects their speed and motion precision. In this paper, a friction isolator (FI) is used as a simple and robust solution to mitigate the undesirable effects of pre-motion friction in precision motion stages. For the first time, a  theoretical study is carried out to understand the dynamic phenomena associated with using a friction isolator on a motion stage. Theoretical analysis and numerical simulation are conducted to examine the dynamical effects of FI  on a PID-controlled motion stage under LuGre friction dynamics. The influence of FI on the response and stability of the system is examined through theoretical and numerical analysis. Parametric analysis is also carried out to study the effects of FI and friction parameters on the eigenvalue and stability characteristics. The numerical results validate the theoretical findings and demonstrate several other interesting nonlinear phenomena associated with the introduction of FI. This motivates deeper nonlinear dynamical analyses of FI for precision motion control. 
\end{abstract}

\blfootnote{The latest version of this manuscript is published in Journal of Vibration and Control. (DOI: 10.1177/1077546321999510)}
\blfootnote{$^{1}$ Department of Mechanical Engineering, Virginia Tech, Blacksburg, VA, 24061, USA}
\blfootnote{$^{2}$ Department of Mechanical Engineering, University of Michigan, Ann Arbor, MI, 48109, USA}
\blfootnote{$^{3}$ Corresponding Author: obarry@vt.edu}

\input{sec_Intro.tex}
\input{sec_Modeling.tex}
\input{sec_LinearAnalysis.tex}
\input{sec_Result1.tex}

\input{sec_Result2.tex}

\input{sec_Result3.tex}
\input{sec_Conclusion.tex}



\bibliographystyle{unsrtnat}
\bibliography{reference}


\end{document}

%% file: sec_Intro.tex
\section{Introduction and Background}
\label{sec:intro}

Motion stages are used for precision positioning in a wide range of manufacturing and metrology-related processes, such as machining, additive manufacturing, and semi-conductor fabrication \pCite{altintas2011machine}. Mechanical bearings (e.g., sliding and especially rolling bearings) are popular in precision motion stages due to their large motion range, high off-axis stiffness, and cost-effectiveness \pCite{altintas2011machine}. Pre-motion friction is a common problem encountered in mechanical-bearing-based motion stages. The adverse effects in performance caused by pre-motion friction feature large tracking errors, long settling times, and stick-slip phenomena \pCite{futami1990nanometer, armstrong1994survey, al2008characterization, marques2016survey}. In practice, a common servo feedback controller for precision motion stages is the proportional-integral-derivative (PID) controller \pCite{hensen2003friction, kim2011design}. While the PID controller is designed to reduce tracking error, the implementation of feedback controllers in the frictional system may result in self-excited limit cycles known as friction-induced vibrations \pCite{oestreich1996bifurcation, hinrichs1998modelling, van1999approximate, hoffmann2002minimal, hensen2003friction, duffour2004instability, hoffmann2007linear,  nakano2009stick, saha2011delayed, kruse2015influence, pascal2017periodic, niknam2019friction}, which will further afflict the control performance.

In many studies, the compensation of unwanted frictional effects was realized with different controllers. \bluenote{The traditional high-gain PID controller possesses some robustness and can quickly overcome frictional effects, but they may also lead to large overshoots and limit cycles \pCite{armstrong1994survey}. Robust controllers such as $H_\infty$ feedback controllers and disturbance-observer-based controllers can effectively attenuate model uncertainties and disturbances \pCite{kempf1999disturbance, zheng2017design, sariyildiz2019disturbance}. However, these controllers may be limited in mitigating the highly nonlinear and volatile pre-motion friction \pCite{kim2009precise, chong2010practical}. Finally, all controllers are affected by practical conditions including computing power, sampling rate, and noise-sensitivities, making complex controllers (e.g., neural network controllers \pCite{ren2008robust, kim2009precise}) less favorable in the application.}

The friction isolator (FI), also known as the compliant joint, is a mechanical device recently proposed to effectively and robustly mitigate pre-motion friction \pCite{dong2017simple, dong2018experimental}. Unlike the rigid connection (i.e., high stiffness) between mechanical bearings and tables as in the conventional motion stages, FI introduces a lower stiffness between the bearing and the table, which isolates the frictional dynamics and makes the table more compliant in the motion direction. Prior works experimentally demonstrated that with the implementation of FI, the compliant motion stage greatly reduced tracking errors, significantly improved robustness towards friction changes \pCite{dong2017simple}, and notably shortened settling time compared to the conventional motion stage (i.e., without FI) \pCite{dong2018experimental}.

Given the remarkable improvements in positioning precision, quickness, and robustness brought by the FI as observed in experiments \pCite{dong2017simple, dong2018experimental}, it is important to fundamentally understand the beneficial and potentially harmful dynamical effects the FI introduces into the precision motion stage. The conventional motion stage is often modeled as a single-body friction oscillator. The dynamics in single-body friction oscillators has been extensively investigated in these works, which feature analyses of stability \pCite{oestreich1996bifurcation, hinrichs1998modelling, hoffmann2007linear, saha2011delayed, kruse2015influence}, mode coupling \pCite{hoffmann2002minimal, hoffmann2007linear}, nonlinear behavior \pCite{oestreich1996bifurcation, feeny1997phase, hinrichs1998modelling, van1999approximate, saha2011delayed, kruse2015influence}, and bifurcation \pCite{oestreich1996bifurcation, hinrichs1998modelling, van1999approximate, di2003sliding,  saha2011delayed}. Implementing FI will introduce new inertia, stiffness, and damping elements into the system. Studies have also been conducted on the friction-induced vibration of multibody systems \pCite{galvanetto1999non, duffour2004instability, nakano2009stick, pascal2017periodic, niknam2019friction}. However, these works either have not considered pre-motion frictional dynamics, or they have adopted models that do not fit with the FI-equipped motion stage system. Furthermore, while the effect of the integral controller on the frictional dynamics of single-body systems has been investigated \pCite{hensen2003friction, bisoffi2017global}, none of the studies has explored the frictional behavior of a PID controlled internally coupled multibody slider system. These problems are investigated for the first time in this paper. The aim is to understand the frictional dynamics of a PID-controlled motion stage system with and without FI under matching parameter conditions. As an extension of our conference paper \pCite{dong2019friction}, the dynamic models are established for PID-controlled motion stages both with and without FI coupled with the LuGre friction model, whose system parameters are experimentally obtained. Numerical simulations are carried out to validate the theoretical analyses.  Parametric studies are conducted to understated the role of PID control gains, friction parameters, and FI design parameters on the stability of the motion stage.

The remaining contents of the paper are organized as follows: 
we first establish and discuss the dynamical models of the PID-controlled stages with and without FI under LuGre friction. Next, the linear analysis of the models is carried out by studying the properties of the state Jacobian matrices. We then validate the findings, and parametrically study the effect of LuGre friction and FI on the performance and stability of PID-controlled motion stages through both analytical and numerical approaches; Finally, we summarized our findings and propose future works in the conclusion section.

%% file: sec_Modeling.tex
\section{Dynamical Modeling of Motion Stages with Friction}
\label{sec:modeling}

\begin{figure}
\begin{center}
\includegraphics[width = 0.6\linewidth]{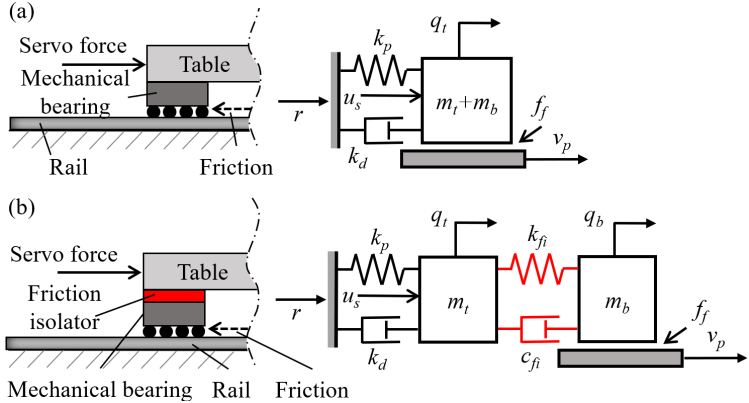}
\caption{Schematics of servo-controlled motion stage under friction - (a) System $\alpha$: the conventional motion stage without FI; and (b) System $\beta$: the compliant motion stage with FI.}
\label{fig:pic_Schematic_Motion_Stage} 
\end{center}
\end{figure}

This section introduces the dynamical modeling of the servo-controlled motion stages systems both with and without FI. For simplicity of labeling, the motion stage without FI is referred to as System $\alpha$, and the motion stage with FI is referred to as System $\beta$ hereinafter. 

Figure \ref{fig:pic_Schematic_Motion_Stage}(a) shows the schematics of a servo-controlled conventional motion stage without FI (System $\alpha$). The mass of the moving table rigidly connected to the bearing is $m_t+m_b$, where $m_t>0$ and $m_b>0$ are the table mass and the bearing mass, respectively. The table is coupled with the reference command $r$ by stiffness $k_p$ and damping $k_d$, which are respectively equivalent to the proportional and derivative gains of the servo feedback controller that regulate the table position $q_t$.  A supplemental control force $u_s$ is added to account for additional servo forces (e.g., feedforward force, integral controller) that may be applied to the table. The friction force $f_f$ is directly applied to the table via a moving platform at velocity $v_p$ .

The compliant motion stage with FI (System $\beta$) is modeled as a system of two coupled oscillators as depicted in Fig. \ref{fig:pic_Schematic_Motion_Stage}(b). While System $\beta$ is largely similar to System $\alpha$, the newly introduced FI connects the table and the bearing via the stiffness $k_{\FI}$ and damping $c_{\FI}$, resulting in the internal coupling between the two bodies. The positions of the table and the bearing are $q_t$ and $q_b$, respectively. The friction force $f_f$ acting on the bearing is not directly applied to the table.

\subsection{The LuGre friction model}
A variety of friction models have been proposed in the past decades. In this study, the LuGre friction model \pCite{de1995new,hoffmann2007linear} is adopted, which incorporates viscous friction, pre-motion friction (i.e., pre-sliding/pre-rolling), and hysteresis behaviors. The LuGre model introduces an internal state $z$, which is used to represent the average deflection of the contact bristles between two surfaces at the friction interface. The dynamics of $z$ is given by
\begin{equation}
    \dot{z}= v - a_z(v) z
\end{equation}
with
\begin{equation}
    a_{z}(v)=\frac{|v|}{g(v)}; 
     \quad  g(v) =  \frac{f_{C} + (f_{S} - f_{C}) e ^{- ({v}/{v_{s}})^2}}{\sigma_{0}}
\end{equation}
where $v$ is the relative velocity between two moving surfaces, $f_C$ is the Coulomb friction, $f_S$ is the static friction, $v_s$ is the Stribeck velocity threshold, and $\sigma_0$ is the initial contact stiffness of the bristle. The modeled friction force of the LuGre model is then calculated as
\begin{equation}
    f_f = {\sigma_0} z + {\sigma_1} \dot{z} + {\sigma_2} v
\end{equation}
where $\sigma_1$ is the micro-damping of the bristle, and $\sigma_2$ accounts for the macroscopic viscous friction. 

Observe that the dynamics of $z$ is only affected by $v$. Hence, the equilibrium points of $z$ can be reached only when
\begin{equation}
    \text{(1): } v=0 \quad \text{or} \quad \text{(2): } z= v/a_z(v)
    \label{eq:LuGreEquilibrium}
\end{equation}
Equilibrium (1) is known as the sticking equilibrium, and equilibrium (2) is referred to as the slipping equilibrium. The fixed points of any dynamic systems that involve the LuGre friction have to satisfy either of these two conditions. It should be noted that the dynamical model is a switched system at $v=0$ due the existence of $\text{sgn}(v)$ and $|v|$.

\subsection{Connections between servo-controlled stage and self-excited friction oscillator}
As discussed before, the models shown in Fig.  \ref{fig:pic_Schematic_Motion_Stage} are often used to study the dynamics of a servo-controlled motion stage during trajectory tracking application, assuming  $r \neq 0$ and $v_p = 0$. In this case, the moving table of the stage (i.e., $q_t$) is controlled to follow a time-varying reference signal $r(t)$. The resulting tracking error can be obtained as 
\begin{equation}
    \epsilon = q_t - r(t);
\end{equation}
When the industrial-standard linear PID controller is implemented, the feedback control force can written as 
\begin{equation}
    u_b = - \epsilon_i -  k_p \epsilon - k_d \dot{\epsilon}
    \label{eq:linearPIDStructure}
\end{equation}
where $\epsilon_i = k_i \int \epsilon dt$; and $k_p$, $k_i$, $k_d \geq 0$ are respectively the proportional, integral, and derivative gains. If we denote $y_i=k_i \int q_t dt$ as the additional state brought about by the integral action, the full states of the systems are defined as
\begin{subequations}
\begin{gather}
    \mathbf{x}_{\alpha}=\bbmat q_t & \dot{q}_t & z & y_{i} \ebmat^\tp
    \\
    \mathbf{x}_\beta = \bbmat q_t & q_b & \dot{q}_t & \dot{q}_b & z & y_{i} \ebmat^\tp
\end{gather}
\end{subequations}
and the state equations of the systems with and without FI are given by
\begin{equation}
    \dot{\mathbf{x}}_\alpha=\bbmat
    \dot{q}_t \\ (-f_f+u)/m_\alpha \\ v-a_z(v) \\ k_i q_t
    \ebmat
    ; \
    \dot{\mathbf{x}}_\beta=\bbmat
    \dot{q}_t \\ \dot{q}_b \\ (f_{\FI}+u)/(m_t) \\ (-f_{\FI}-f_f)/(m_b) \\ v-a_z(v) \\ k_i q_t
    \ebmat
    \label{eq:originalSS}
\end{equation}
where 
\begin{equation}
    f_{\FI} = k_{\FI}( q_b - q_t ) + c_{\FI} ( \dot{q}_b - \dot{q}_t )
\end{equation}
is the coupling force of the friction isolator, $m_\alpha=m_t+m_b$ is the total mass of the system, and $u$ is the additional control input (e.g., feedforward action). The above state equations can be re-arranged to obtain the error dynamics as
\begin{subequations}
\begin{gather}
    \dot{\mathbf{x}}_\alpha^\star=\bbmat
    \dot{\epsilon} \\ (-f_f+u)/m_\alpha - \Ddot{r} \\ v-a_z(v) \\ k_i \epsilon
    \ebmat
    \\
    \dot{\mathbf{x}}_\beta^\star=\bbmat
    \dot{\epsilon} \\ \dot{\epsilon}_b \\ (f_{\FI}+u)/m_t - \Ddot{r} \\ (-f_{\FI}-f_f)/m_b - \Ddot{r} \\ v-a_z(v) \\ k_i \epsilon
    \ebmat
\end{gather}
    \label{eq:alternativeSS}
\end{subequations}
where 
\begin{subequations}
\begin{gather}
    \mathbf{x}_{\alpha}^\star=\bbmat \epsilon & \dot{\epsilon} & z & \epsilon_i \ebmat^\tp 
     \\ 
    \mathbf{x}_\beta ^\star= \bbmat \epsilon & \epsilon_b & \dot{\epsilon} & \dot{\epsilon}_b & z & \epsilon_i \ebmat^\tp 
     \\ 
    \epsilon_b = q_b - r
\end{gather}
\end{subequations}

In the literature of friction-induced vibrations under self-excitation \pCite{hinrichs1998modelling,hoffmann2007linear,li2016nonlinear,saha2016modified}, friction is often introduced by fixing the reference and prescribing the platform (e.g., belt, conveyor) with constant velocity motion (i.e., $r = 0$, $v_p \neq 0$, $\dot{v}_p = 0$), as shown in Fig. \ref{fig:pic_Schematic_Motion_Stage}. Therefore, the relative velocity $v$ between the frictional interfaces can be written as
\begin{equation}
    v_\alpha= \dot{q}_t-v_p; \quad v_\beta=\dot{q}_b-v_p
    \label{eq:friction oscillator}
\end{equation}
In the case of servo-controlled motion stage (i.e., $r \neq 0$, $v_p = 0$), the relative velocities are obtained as,
\begin{equation}
    v_\alpha=\dot{q_t} =  \dot{\epsilon} + v_r; \quad v_\beta=\dot{q_b} = \dot{\epsilon}_b+v_r
    \label{eq:motion_control}
\end{equation}

Note that by substituting $r = 0$ in Eq.(\ref{eq:friction oscillator}) (that is, $\dot{q}_t = \dot{\epsilon}$, $\dot{q}_b = \dot{\epsilon}_b$) and setting $v_r = \dot{r} = - v_p$ (i.e., reference trajectory is constant velocity motion) in Eq.(\ref{eq:motion_control}), the relative velocities of these two cases become identical. This indicates that the dynamical response and stability of self-excited friction oscillator and servo-controlled motion stage are equivalent. Even when $\ddot{r} \neq 0$, the equivalence can be acquired simply by designing the controller $u$ as
\begin{equation}
    u = u_f + u_b; \quad u_f= m_\alpha \ddot{r}
    \label{eq:fullController}
\end{equation}
where $u_f$ is the feedforward controller that provide the acceleration. Therefore, in the rest of this paper, we will focus on the representation pertaining to the servo-controlled motion stage.

%% file: sec_LinearAnalysis.tex
\section{Linear Stability Analysis}
\label{sec:linearStability}
The introduction of FI to the servo-controlled motion stage may pose challenges to the stability of the system. In this section, the effects of FI and friction parameters on the feedback controller design of the stage are investigated using linear stability analysis. 

\subsection{Equilibrium points and state Jacobian matrices}
As part of the stability analysis, the calculation of the state equilibrium points may vary with the state representations of the system. For Systems $\alpha$ and $\beta$, the equilibrium is studied with respect to the error dynamics (Eq.(\ref{eq:alternativeSS})) and controller dynamics (Eq.(\ref{eq:linearPIDStructure})). As mentioned in the previous section, the equilibrium of the system with the LuGre friction can be reached only when either of the two conditions in Eq.(\ref{eq:LuGreEquilibrium}) is satisfied. When the PD controller is applied ($k_i=0$), the stick equilibrium points (at $v_r=0$) are calculated as
\begin{subequations}
    \begin{gather}
        \mathbf{x}_{\alpha,0}^\star = \bbmat \epsilon_{0} & 0 & -k_p\epsilon_{0}/\sigma_0 & 0\ebmat^\tp
        \\
        \mathbf{x}_{\beta,0}^\star = \bbmat \epsilon_{0} &  ( k_t \epsilon_{0}- m_b \Ddot{r} ) / k_{\FI} & 0 & 0 & -k_p \epsilon_{0}/\sigma_0 & 0\ebmat^\tp
    \end{gather}
\end{subequations}
where $k_t = k_{p} + k_{\FI} $; and $\epsilon_0 \in \mathbb{R}$ is the steady state position error. When $k_i \neq 0$, the sticking equilibrium points become
\begin{subequations}
    \label{eq:stickingEquilibriumPID}
    \begin{gather}
        \mathbf{x}_{\alpha,0}^\star = \bbmat 0 & 0 & -\epsilon_{i,0}/(\sigma_0) & \epsilon_{i,0}\ebmat^\tp
        \\
        \mathbf{x}_{\beta,0}^\star = \bbmat 0 &  (\epsilon_{i,0} - m_b \Ddot{r} ) / k_{\FI} & 0 & 0 & -\epsilon_{i,0}/\sigma_0  & \epsilon_{i,0} \ebmat^\tp
    \end{gather}
\end{subequations}
where $\epsilon_{i,0} \in \mathbb{R}$ is the integral error that balances the unmodeled system dynamics. 

The slipping equilibrium occurs when $v_r \neq 0$. In the absence of integral action (i.e., PD control), the equilibrium points are obtained as
\begin{equation}
        \mathbf{x}_{\alpha,0}^\star = \bbmat -f_{f,0} /k_p \\ 0 \\ h(v_r) \\ 0\ebmat
        ; \quad
        \mathbf{x}_{\beta,0}^\star = \bbmat -f_{f,0} /k_p \\  \epsilon_{b,0,pd} \\ 0 \\ 0 \\ h(v_r) \\ 0\ebmat
\end{equation}
where
\begin{gather*}
    h(v)= v/a_z(v) = \text{sgn}(v)g(v); \ \  f_{f,0}=\sigma_0 h(v_r) +\sigma_2 v_r \\ 
    \epsilon_{b,0,pd} = - ( k_t f_{f,0} + k_p m_b \Ddot{r}) / ( k_p k_{\FI}) 
\end{gather*}
Similarly, the slipping equilibrium points in the presence of PID controllers can be calculated as
\begin{equation}
    \mathbf{x}_{\alpha,0}^\star = \bbmat 0 \\ 0 \\  h(v_r)\\ -f_{f,0}\ebmat
    ; \
    \mathbf{x}_{\beta,0}^\star = \bbmat 0 \\ -( f_{f,0} + m_b \Ddot{r} ) / k_{\FI} \\ 0 \\ 0 \\ h(v_r) \\ - f_{f,0}\ebmat
\end{equation}

Linear stability analysis is carried out by examining the Hurwitz properties of the state Jacobian matrix, which is obtained by linearizing the system around the equilibrium points \pCite{ogata2002modern,khalil2002nonlinear}. \bluenote{To make sure that the steady-state solution of a nonlinear system locally converges to a fixed-point, it is necessary for the linearized system at the fixed point to be stable.} The sticking equilibrium assumes $v_r = 0$, which is not relevant to the scope of this paper. Therefore, the stability analysis is conducted at the slipping equilibrium points. For System $\alpha$ and $k_i=0$, the state Jacobian matrix is
\begin{equation}
    \mathbf{A}_\alpha 
    = \bbmat
        0 & 1 & 0 \\
        -k_p/m_a & a_{\alpha,[2,2]} & a_{\alpha,v,z}\\
        0 & a_{\alpha,z,v} & -a_z(v) \\
    \ebmat
\end{equation}
with
\begin{subequations}
\begin{gather}
    a_{\alpha,[2,2]} = -(k_d + \sigma_1 + \sigma_2 - \sigma_1 z (\partial a_z(v)/ \partial v))/m_a
    \\
    a_{\alpha,v,z} = -(\sigma_0 - \sigma_1 a_z(v))/m_a
    \\
    a_{\alpha,z,v}= 1 - z (\partial a_z(v)/ \partial v)
\end{gather}
    \label{eq:Jacobian_elements_PD}
\end{subequations}
where
\begin{equation}
    \frac{\partial a_z(v)}{\partial v} = \frac{ \text{sgn}(v) [g(v) v_{s}^2 + 2v^2 (g(v) - f_{C}/\sigma_0) ]}{g^2(v)v_{s}^2}
\end{equation}
Since $z=h(v)$ at the equilibrium, Eq.(\ref{eq:Jacobian_elements_PD}) can be further simplified as
\begin{subequations}
\begin{gather}
     a_{\alpha,[2,2]} = - (k_d + \sigma_2 - \sigma_1 \rho_f (v) {v^2}/{v_s^2})
     \\
     a_{\alpha,z,v} = - \rho_f (v) {v^2}/{v_s^2}
\end{gather}
\end{subequations}
where
\begin{equation}
    \rho_f (v) =2- \frac{2f_C}{f_C+(f_S-f_C)e^{-(v/v_s)^2}}
    \label{eq:couplingTerm}
\end{equation}
is the ratio bounded by $(0,2(f_S-f_C)/f_S]$. Similarly, for System $\beta$ and $k_i=0$, the Jacobian can be calculated as
\begin{equation}
    \mathbf{A}_\beta 
        = \bbmat
        0 & 0 & 1 & 0 & 0 \\
        0 & 0 & 0 & 1 & 0 \\
        a_{\beta,[3,1]}  & k_{\FI}/m_t & a_{\beta,[3,3]} & c_{\FI}/m_t & 0\\
        k_{\FI}/m_b & -k_{\FI}/m_b & c_{\FI}/m_b & a_{\beta,[4,4]} & a_{\beta,v,z}\\
        0 & 0 & 0 & a_{\beta,z,v} & -a_z(v)\\
    \ebmat
\end{equation}
where
\begin{subequations}
\begin{gather}
    a_{\beta,[3,1]}=-(k_{\FI}+k_p)/m_t
    \\
    a_{\beta,[3,3]}=-(c_{\FI}+k_d)/m_t
    \\
    a_{\beta,[4,4]} = -(c_{\FI} +\sigma_2  - \sigma_1 \rho_f (v) {v^2}/{v_s^2})/m_b;
    \\
    a_{\beta,v,z} = -(\sigma_0 - \sigma_1 a_z(v))/m_b;
    \\
    a_{\beta,z,v}= - \rho_f (v) {v^2}/{v_s^2}
\end{gather}
\end{subequations}
Notice that the sign of $v$ does not affect the values of state Jacobian matrices, thus confirming the symmetry property of the system. These matrices are only dependent on the states $\dot{\epsilon}$ (for System $\alpha$), $\dot{\epsilon}_b$ (for System $\beta$), and $z$. The resulting values of Jacobian matrices at the slipping equilibrium points are obtained by setting $v=v_r$. For the PID cases, the state Jacobian matrices can be defined as
\begin{equation}
    \mathbf{A}_{\alpha,i} = \bbmat
    0 & \bbmat k_i & 0 & 0 \ebmat\\
    \bbmat  0 -1/m_a & 0\ebmat^\tp & \mathbf{A}_\alpha 
    \ebmat
\end{equation}
and
\begin{equation}
    \mathbf{A}_{\beta,i} = \bbmat
    0 &  \bbmat k_i & 0 & 0 & 0 & 0  \ebmat \\
    \bbmat 0 & 0 & -1/m_t & 0 & 0\ebmat^\tp & \mathbf{A}_\beta
    \ebmat
\end{equation}
which are obtained by re-arranging the sequence of the states (i.e., moving $\epsilon_i$ to the first state).

Finally, it is helpful to convert the dimensional Jacobian matrix to a non-dimensional form such that the eigenvalues are scaled for easier comparison. The general procedure is to select a principal natural frequency $\omega_n$ and use the corresponding non-dimensional time $t_n =\omega_n t$. In this study, the principal natural frequency for the two systems are selected as
\begin{equation}
    \omega_{n,\alpha} = \sqrt{k_p/m_\alpha};\quad \omega_{n,\beta} = \sqrt{k_p/m_t}
\end{equation}
The non-dimensional state Jacobian matrices of the systems can then be obtained as 
\begin{subequations}
\begin{gather}
    \mathbf{A}_{n,\alpha,i} = \mathbf{\Omega}_{n,\alpha,1} \mathbf{A}_{\alpha,i} \mathbf{\Omega}_{n,\alpha,2}
    \\
    \mathbf{A}_{n,\beta,i} = \mathbf{\Omega}_{n,\beta,1} \mathbf{A}_{\beta,i} \mathbf{\Omega}_{n,\beta,2}
\end{gather}
\end{subequations}
where
\begin{subequations}
\begin{gather}
    \mathbf{\Omega}_{n,\alpha,1}= \text{diag}(\bbmat 1 & \omega_{n,\alpha}^{-1} & \omega_{n,\alpha}^{-2} & \omega_{n,\alpha}^{-1}\ebmat)
    \\
    \mathbf{\Omega}_{n,\alpha,2}= \text{diag}(\bbmat \omega_{n,\alpha}^{-1} & 1 & \omega_{n,\alpha} & 1\ebmat)
    \\
    \mathbf{\Omega}_{n,\beta,1}= \text{diag}(\bbmat 1 & \omega_{n,\beta}^{-1} & \omega_{n,\beta}^{-1} & \omega_{n,\beta}^{-2}  & \omega_{n,\beta}^{-2} & \omega_{n,\beta}^{-1}\ebmat)
    \\
    \mathbf{\Omega}_{n,\beta,2}= \text{diag}(\bbmat \omega_{n,\beta}^{-1} & 1 & 1 & \omega_{n,\beta} & \omega_{n,\beta} & 1\ebmat)
\end{gather}
\end{subequations}
Note that the stability implied by the non-dimensional Jacobian matrices are identical to those from the dimensional matrices.

\subsection{Properties of state Jacobian matrices}

\label{subsec:JacobianProperties}
To guarantee the stability of the linearized system, the state matrix has to be Hurwitz, i.e., all eigenvalues have negative real parts. This can be evaluated by directly calculating the eigenvalues or applying the Routh-Hurwitz criterion \pCite{ogata2002modern} on the characteristic equations calculated from the state Jacobian matrices. For example, the characteristic equation of $\mathbf{A}_\alpha$ can be calculated as
\begin{equation}
    s^3 + b_1 s^2 + b_2 s + b_3 = 0
\end{equation}
where the roots of this characteristic equation are the eigenvalues, and 
\begin{subequations}
\begin{gather}
    b_1=a_z(v) + [k_d + \sigma_2 - \sigma_1 \rho_f(v) v^2/ v_s^2]/m_\alpha 
    \\
    b_2= [k_p + k_d a_z(v) + \sigma_2 a_z(v)- \sigma_0 \rho_f(v) v^2/ v_s^2] / m_\alpha
    \\
    b_3= a_z(v) k_p/ m_\alpha
\end{gather}
    \label{eq:routhParam}
\end{subequations}
{If $b_1 b_2 > b_3$, the system is linearly stable. While a complete symbolical evaluation of the Hurwitz property is very difficult due to the complexity of the system}, for System $\alpha$ under PD control, it is observed that when $\sigma_0 \rightarrow \infty$, other coefficients have trivial effects on the stability of the system. This leads  to $ b_1 \approx a_z(v)$ and the simplified stability condition can be obtained as
\begin{equation}
    b_1 b_2 - b_3 \approx \frac{a_z^2(v)}{m_\alpha^2}(k_d+\sigma_2 - \frac{|v| (f_S-f_C)e^{-(\frac{v}{v_s})^2}}{v_s^2}) >0
    \label{eq:extremeSigma}
\end{equation}
Note that $\sigma_1$ is not presented in the above expression. This indicates that micro-damping does not affect system stability when $\sigma_0 \rightarrow \infty$. Also, smaller $v_s$ requires larger $k_d$ or $\sigma_2$ to stabilize the system, especially when $v$ is close to $v_s$. In addition, for a fixed $v_s$, the lower bounds of $k_d$ and $\sigma_2$ that stabilize the system reach maximum values when $v=v_s/\sqrt{2}$; this is obtained by taking the derivative of ${|v|(f_S-f_C)e^{-({v}/{v_s})^2}}/{v_s^2}$ with respect to $v$. 

Several other interesting properties of the systems can also be obtained by examining the structure of the Jacobian matrices. Notice that all the  matrices can be decomposed into the following structure
\begin{equation}
    \mathbf{A}=\bbmat \mathbf{A}_{M} & \bbmat 0 & \cdots & 0 & a_{v,z} \ebmat \\
    \bbmat 0 & \cdots & 0 & a_{z,v} \ebmat^\tp & -a_z(v)
    \ebmat
\end{equation}
where $\mathbf{A}_{M}$ is the submatrix corresponding to the states from the multibody system. For both systems, it can be observed that $\rho_f(v) \rightarrow 0$ when $|v| \gg v_s$. As a result, when the relative velocity at friction surface is significantly larger than the Stribeck velocity threshold, $a_{z,v} \rightarrow 0$ and the eigenvalues of $\mathbf{A}$  consist of eigenvalues of $\mathbf{A}_{M}$ and $-a_z(v)$; the latter is negative by default. This indicates that the effects of friction dynamics and multibody dynamics on the system stability are decoupled when the  velocity is large. The same conclusion can be drawn when $v_s \gg |v| >0$ or $f_S \approx f_C$, which all lead to $a_{z,v} \rightarrow 0 $. Based on this property, the analysis scope can be reasonably focused on the low speed range (i.e., when $v_{r}$ is close to $v_{s}$). Similarly, the structures of $\mathbf{A}_{i}$ indicate that, when $k_i$ is small, the eigenvalues from $\mathbf{A}$ will be unaffected and carried over to $\mathbf{A}_{i}$. This property makes it convenient to distinguish the eigenvalue introduced by the additional integral state $\epsilon_i$. 

%% file: sec_Result1.tex
\section{Results and Discussion}
\label{sec:simulation}

In this section, we first validate the theoretical observations using numerical simulation, and then we examine the effects of LuGre friction and FI on the performance and stability of PID controlled motion stages. {The default design and friction parameters obtained from the prototype in previous experimental studies \pCite{dong2017simple, dong2018experimental, dong2019friction} are:
\begin{subequations}
\begin{gather}
    m_t = 1 \text{ kg}; \ m_b = 0.5 \text{ kg}; \ f_S = 6.5 \text{ N}
    \\
     \quad f_C = 5.1 \text{ N}; \ k_{\FI} = 4 \times 10^4 \text{ N/m}; \ c_{\FI} = 2 \text{ Ns/m}
    \\
    v_s = 16.7 \text{mm/s}; \ \sigma_0 = 2.2 \times 10^6 \text{N/m}
    \\
    \sigma_1 = 237 \text{Ns/m}; \  \sigma_2 = 14.2 \text{Ns/m}
\end{gather}
   \label{eq:standardParam}
\end{subequations}
The default PID controller gains are:
\begin{subequations}
\begin{gather}
     k_p= 2\times 10^4 \text{ N/m}; \quad k_i= 1 \times 10^6 \text{ N/ms}
     \\ 
     k_d= 2\times 10^2 \text{ Ns/m}; 
\end{gather}
   \label{eq:standardGain}
\end{subequations}
The default reference velocity is chosen as ${v}_r = 10 \text{mm/s}$. The initial conditions of the numerical simulations are selected as $\mathbf{x}_\alpha^\star=[0,-v_r,0,0]^\tp$ and $\mathbf{x}_\beta^\star=[0,0,-v_r,-v_r,0,0]^\tp$.}

\subsection{Numerical Validation of the Theoretical Results}

To validate the results from the linear stability analysis, numerical simulations are first carried out using the built-in ODE solver \texttt{ode45} in MATLAB with the nonlinear system equations. 
Constant velocity motion is used as the reference trajectory, i.e., $r=v_r t$. The parameters used in the numerical validation are from Eq.(\ref{eq:standardParam}, \ref{eq:standardGain}).

\begin{figure}[!b]
    \centering
    \includegraphics[width = 0.6\linewidth]{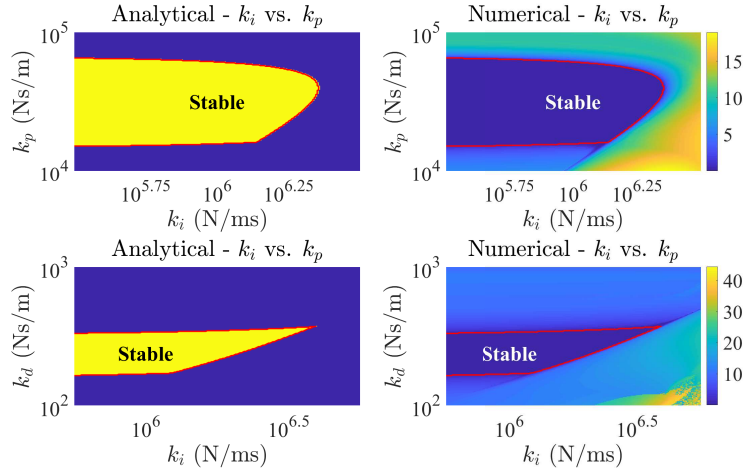}
    \caption{The analytical stability region of PID control gains and its corresponding steady state oscillation error amplitude map acquired through the numerical simulation. \bluenote{For analytical results, the red contours are stability boundaries; for numerical results, the red contours enclose the steady state solutions that converge to zeros (i.e., fixed points).}
    }
    \label{fig:NumericalValidation}
\end{figure}

\begin{figure}[!t]
    \centering
    \includegraphics[width = 0.6\linewidth]{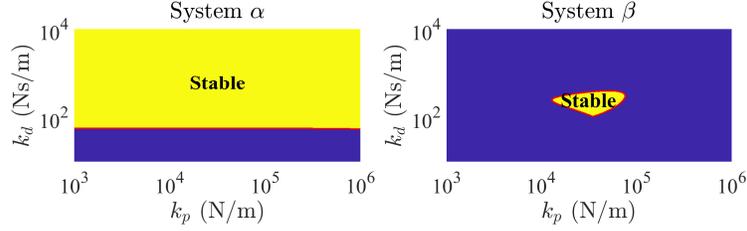}
    \caption{Region of stable PD gain for System $\alpha$ and $\beta$ at $v_r=10$ mm/s, where the red contour is the boundary of the stable domain.}
    \label{fig:stabilityAreaPD2D}
\end{figure}

It is observed that the numerical stability analysis closely matches with the one obtained analytically. {An example carried out on System $\beta$ is demonstrated in Fig.  \ref{fig:NumericalValidation}, which compares the linear stability regions with the steady state oscillation amplitudes respectively within the $k_i$-$k_p$ domain and the $k_i$-$k_d$ domain. The stability regions are calculated by evaluating the Hurwitzness of the state matrix $\mathbf{A}_{\beta,i}$ through eigenvalues. The steady state oscillation amplitudes are calculated by simulating the ODE of System $\beta$ until the transient responses are phased out. 

The analytical results yield a boundary that separates the stable and unstable regions. Although the numerical results do not provide a direct indication of stability, \bluenote{the steady state solutions of $\mathbf{x}^\star_\beta$ that converge to zeros (i.e., fixed-points) can be separated from the ones with non-zero amplitudes (i.e., oscillations).} The separation contours match excellently with the stability boundaries calculated analytically. This is an indication that the stability of the nonlinear system can be reliably evaluated through the analytical approach from the linear analysis, which leads to the following parametric study on the stability of PID-controlled motion stages.}

%% file: sec_Result2.tex
\subsection{{Effect of Friction on System Stability}}
\label{subsec:ctrlParamStudy}

\begin{figure}[!b]
    \centering
    \includegraphics[width = 0.6\linewidth]{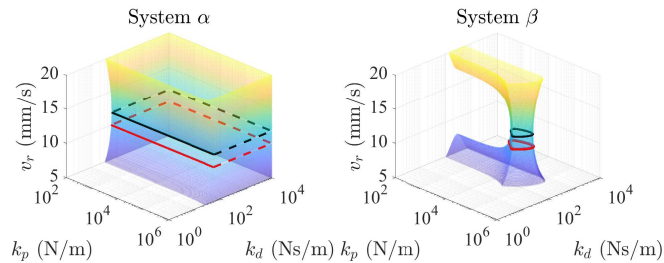}
    \caption{The evolution of the boundary of the stable PD gain domain (enclosed by the boundary) with respect to $v_r$ for both systems, where the red contour is the boundary at $v_r=10$ mm/s (also observable in Fig. \ref{fig:stabilityAreaPD2D}), and the black contour is the boundary at $v_r=v_s/\sqrt{2} \approx 11.8$ mm/s. The dash-line edges indicate the openings of the boundary  (i.e., the system does not become unstable with very high $k_d$ gains or $k_p$ gains).}
    \label{fig:stabilityAreaPD3DAlongv0}
\end{figure}

\begin{figure}[!t]
    \centering
    \includegraphics[width = 0.6\linewidth]{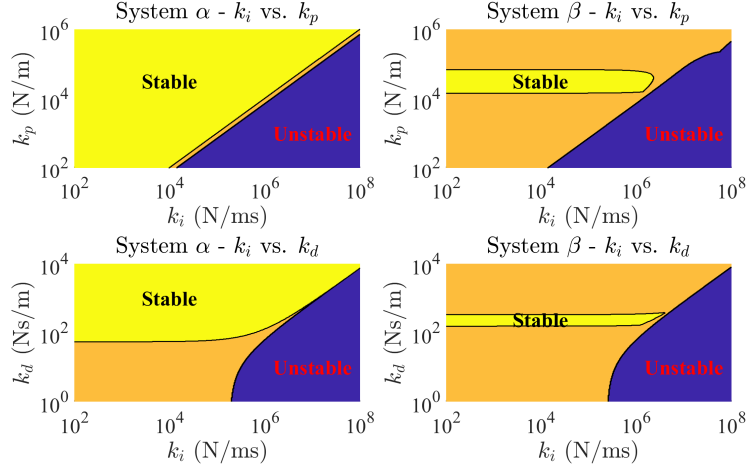}
    \caption{Stability region of System $\alpha$ and $\beta$ under PID control, where: yellow areas indicate stability for both $v_r=10$ mm/s and $|v_r| \gg v_s$ (i.e. friction dynamics is isolated) cases; orange areas indicate instability for $v_r =10$ mm/s but stability when $|v_r| \gg v_s$; and blue areas indicate instability for both $v_r =10$ mm/s and $|v_r| \gg v_s$.}
    \label{fig:PIDstabilityComparison}
\end{figure}

\begin{figure}[!b]
    \centering
    \includegraphics[width = 0.6\linewidth]{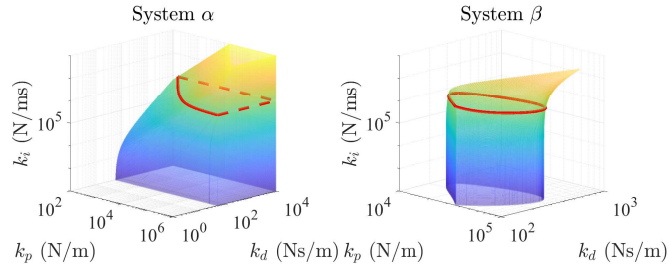}
    \caption{The evolution of the boundary of the stable $k_p$-$k_d$ domain (enclosed by the boundary) with respect to $k_i$ for both systems, where the red contour is the boundary at the default $k_i=1 \times 10^6$ N/ms. The dash-line edges indicate the openings of the boundary  (i.e., the system does not become unstable with very high $k_d$ gains or $k_p$ gains)}
    \label{fig:stabilityAreaPD3DAlongki}
\end{figure}

In the servo-controlled motion stage, the PID controller is designed to first stabilize the system. Therefore, the effects of friction on the stability of PID-controlled motion stages are investigated with a focus on how the ranges of stable control gains are affected. As discussed in the modeling section, the LuGre friction force $f_f$ is determined by a total of seven parameters - $v$, $v_s$, $f_S$, $f_C$, $\sigma_0$, $\sigma_1$, and $\sigma_2$, among which only $v$ is state dependent. Based on the fact that $v=v_r$ at the equilibrium points, the ranges of stable PD control gains (i.e., $k_i = 0$) at the default reference velocity $v_r=10$ mm/s are compared between System $\alpha$ and $\beta$ - see Fig. \ref{fig:stabilityAreaPD2D}. It can be observed that a decrease in the stable gain region occurs when FI is applied. Similar observation is also shown in Fig. \ref{fig:stabilityAreaPD3DAlongv0}, where the variation of the stable gain boundary with different $v_r$ is presented. In general, it is much harder to tune the PD gains that can stabilize System $\beta$ to the equilibrium point. For System $\alpha$, the boundary can be predicted by Eq.(\ref{eq:extremeSigma}). The results suggest that $k_d$ needs to be larger than a certain value to overcome the destabilizing effect of friction. It is also observed that at the critical velocity $v_r=v_s/\sqrt{2} \approx 11.8$ mm/s, System $\alpha$ has the smallest stable region. {This matches with the finding from Eq.(\ref{eq:extremeSigma}).} Although it is harder to obtain such a value for System $\beta$, Figure \ref{fig:stabilityAreaPD3DAlongv0} shows that the smallest stable gain boundary also appears around $v_r=v_s/\sqrt{2}$. The following analysis is carried out at $v_r$ = 10 mm/s for simplicity. Notice from Fig. \ref{fig:stabilityAreaPD3DAlongv0} that at $v_r$ = 10 mm/s, the stability boundary is close to that from the worst case scenario observed at the critical velocity $v_r=v_s/\sqrt{2}$.

In the presence of integral action, the effects of friction parameters on the ranges of stable control gains are shown in Fig. \ref{fig:PIDstabilityComparison}. In general, an extremely large $k_i$ leads to instability, which is true for both systems in this study. For system $\alpha$, increasing $k_p$ and $k_d$ allows the tuning of a larger stable $k_i$. Note that $k_p(k_d+\sigma_2)>k_i m_\alpha$ is the stability criterion when the LuGre dynamics is decoupled from the rigid-body dynamics (i.e., $|v_r| \gg v_s$ or $v_s \gg |v_r|$). The coupling with friction reduces the stability boundaries in both motion stages (i.e., the orange areas in Fig. \ref{fig:PIDstabilityComparison} are subtracted from the original stable region), which is particularly significant for System $\beta$. Finally, note that the stability boundaries are hardly affected within the range where $k_i$ is small.

Figure \ref{fig:stabilityAreaPD3DAlongki} shows the stable $k_p$-$k_d$ boundaries as the integral gain $k_i$ changes. Observe that the increase of $k_i$ reduces and shifts the stable PD control gain domains for both systems. This indicates that extra care should be taken when using a $k_i$ of large magnitude to quickly overcome the disturbance. For System $\beta$, the stability boundaries of the PD control gains are not much affected by the default value of $k_i$ from Eq.(\ref{eq:standardGain}) when compared with the PD control case (where $k_i=0$).

To summarize, the friction can cause instabilities of PID-controlled motion stages both with and without FI. The observations of the effect of the tracking velocity $v_r$ on the stability boundaries corroborate the analytical findings from the previous section. It is also observed that FI, with the default parameters, can further reduce the ranges of stable PID gains for System $\beta$. This indicates that the choice of FI parameters is critical to the controller design and the stability of the motion stage. 

\subsection{{Effect of Friction Isolator on System Stability}}

An interesting phenomenon observed in the previous subsection is that System $\beta$ experiences unstable-stable-unstable transition as $k_p$ or $k_d$ increases. This is different from System $\alpha$ in which stability is guaranteed once $k_p$ or $k_d$ is larger than a critical value. To further investigate this, the trajectories of the eigenvalues of the systems are evaluated using root locus plots. Based on the default parameters in Eq.(\ref{eq:standardParam}, \ref{eq:standardGain}), it is observed that the following eigenvalues are critical to system stability:
\begin{enumerate}[{(}1{)}]
    \item An eigenvalue $\lambda_i$ (for both systems) introduced by the additional integral state $\epsilon_i$, i.e., $\lambda_i=0$ at $k_i=0$.
    \item A negative real eigenvalue $\lambda_z$ (for both systems) that is introduced by the bristle dynamics with large magnitude.
    \item Eigenvalues $\lambda_t$ and $\lambda_b$ (for System $\beta$ alone) in complex pairs that are introduced by the table and bearing.
\end{enumerate}

\begin{figure}[!t]
    \centering
    \includegraphics[width = 0.6\linewidth]{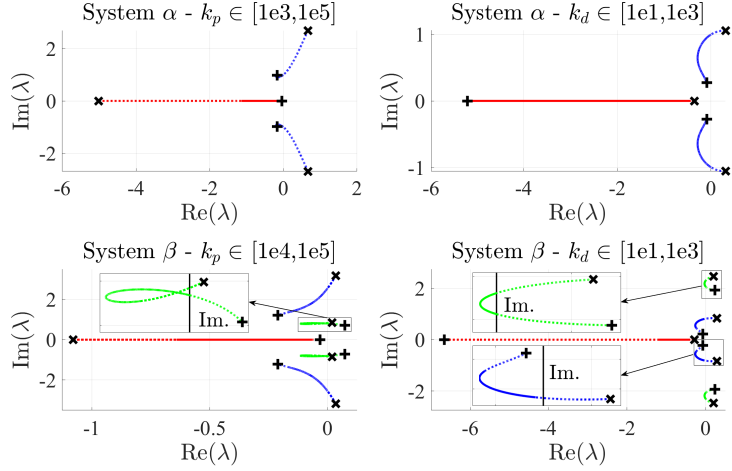}
    \caption{Non-dimensionalized root locus (excluding $\lambda_z$) of System $\alpha$ and $\beta$ with respect to $k_p$ and $k_d$ with $v_r=10$ mm/s and $k_i=1 \times 10^6$  N/ms, where: different colors are used to distinguish eigenvalues; \bluenote{the dotted and solid lines indicate the range of the locus where the systems are unstable (i.e., one or more eigenvalues have positive real parts) and stable (i.e., all eigenvalues have negative real parts)}, respectively; the "x" and "+" markers respectively indicate the beginning and the end of the locus. In the zoomed-in subplots, the black line is the imaginary axis (marked with "Im.").}
    \label{fig:locusAlpha}
\end{figure}

\begin{figure}[!b]
    \centering
    \includegraphics[width = 0.6\linewidth]{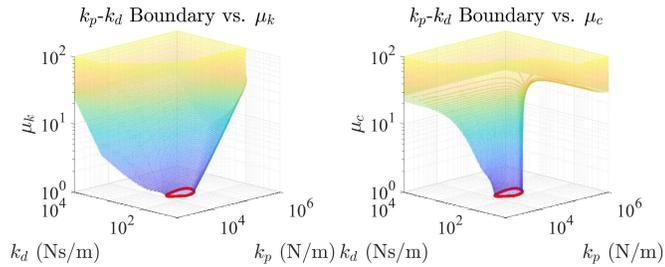}
    \caption{The development of boundary of the stable $k_p$-$k_d$ region (enclosed by the boundary) with respect to $\mu_k$ and $\mu_c$ for System $\beta$ with $v_r=10$ mm/s and $k_i=1 \times 10^6$  N/ms, where the red contour is the boundary with the default parameters.}
    \label{fig:stabilityAreaPD3DAlongfi}
\end{figure}

\begin{figure}[!t]
    \centering
    \includegraphics[width = 0.6\linewidth]{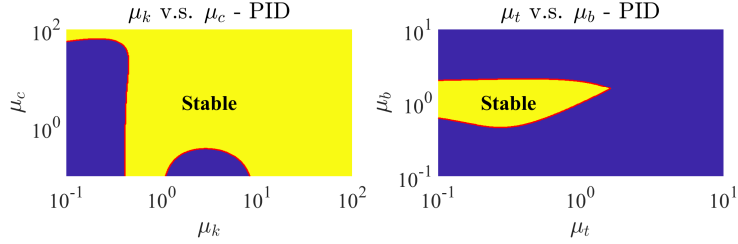}
    \caption{Stability of the system at different design parameters with the default PID gains and $v_r$ defined in Eq.(\ref{eq:standardGain}), where the red contour is the boundary of the stable domain.}
    \label{fig:designParamStability}
\end{figure}

\begin{figure}[!b]
    \centering
    \includegraphics[width = 0.6\linewidth]{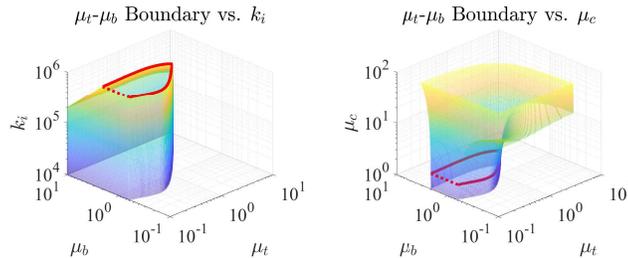}
    \caption{The evolution of the boundary of stable $\mu_t$-$\mu_b$ range (enclosed by the boundary) with respect to $k_i$ and $\mu_c$, where the red contour is the boundary with the default parameters observable from Fig. \ref{fig:designParamStability}. The dot-line edges indicate the openings of boundary (i.e., smaller $\mu_t$ will not lead to instability).}
    \label{fig:massStability}
\end{figure}

The root locus of System $\alpha$ for the abovementioned eigenvalues is presented in Fig. \ref{fig:locusAlpha}. Since $\lambda_z$ has a very large magnitude compared to the other eigenvalues, its trajectories have been excluded from the plot. The figure demonstrates the transition of system stability from unstable to stable as $k_p$ and $k_d$ increase. Note that for System $\beta$, two eigenvalue trajectories (i.e., in blue and red) are very similar to those in the root locus of System $\alpha$; they show the same unstable-stable transition as $k_p$ and $k_d$ increases. However, the additional trajectories of the complex eigenvalue pairs (in green) demonstrates the unstable-stable-unstable transition as they both cross the imaginary axis twice. This explains why the presence of FI shrinks the range of $k_p$ and $k_d$ values that can stabilize the system. {An explanation of the stable-to-unstable transition in System $\beta$ is: when $k_p$ and $k_d$ become very large, the table is rigidly ``constrained" to the tracking reference, resulting in $q_t \approx r$. This leads to System $\beta$ being reduced to a one-body system like System $\alpha$, where $k_{\FI}$ and $c_{\FI}$ play the roles of $k_p$ and $k_d$, respectively. In the current case, $c_{\FI}$ is small, which leads to instability of the system due to the frictional effect.}

Therefore, it is very important to study the effects of FI parameters on system stability. By defining $\mu_k, \ \mu_c \in \mathbb{R}_+$, the FI stiffness and damping can be scaled as $k_{\FI}=\mu_k k_{\FI,0}$ and $c_{\FI}=\mu_c c_{\FI,0}$, where $k_{\FI,0}$ and $c_{\FI,0}$ are the default values from Eq.(\ref{eq:standardParam}). The resulting stable $k_p$-$k_d$ boundary with respect to different scaling coefficients are shown in Fig. \ref{fig:stabilityAreaPD3DAlongfi}. The result indicates that increasing $k_{\FI}$ and $c_{\FI}$ both extend the range of stable PID controller gains. {Note that if $k_{\FI}$ and $c_{\FI}$ are extremely large, the table and bearing will be rigidly connected, which also reduces System $\beta$ to System $\alpha$.}

Previous results show the change of stable PD regions as the FI parameters vary. Alternatively, the influence of design parameters can be visualized by examining their stable combinations with a fixed set of controller parameters. Two new coefficients $\mu_t, \ \mu_b \in \mathbb{R}_+$ are defined so that the mass of the table and bearing can be scaled as $m_t=\mu_t m_{t,0}$ and $m_b=\mu_b m_{b,0}$ respectively. The corresponding stability charts in the $\mu_k-\mu_c$ and $\mu_t-\mu_b$ domains are shown in Fig. \ref{fig:designParamStability} with the default controller gains. Similar to the previous analysis, it is observed that the stability is improved by increasing $k_{\FI}$ and $c_{\FI}$. The results on $\mu_t$ and $\mu_b$  show that the range of stable $m_t$-$m_b$ combinations is quite narrow due to the introduction of FI. This may limit the payload that a motion stage can handle in practice. This problem may be alleviated by increasing the FI damping $c_{\FI}$ or adopting smaller integral gain $k_i$ as shown in Fig. \ref{fig:massStability}. 

In summary, the choice of FI parameters is essential to the stability of the system. There is a trade-off between better isolation performance (in terms of mitigating undesirable effects of pre-motion friction) and improved system stability (in terms of stable regions of PID controller gains) \pCite{dong2017simple,dong2018experimental}. The linear stability analyses from this section can serve as useful guidelines during the design optimizations of FI.

%% file: sec_Result3.tex
\subsection{{Effect of FI on the Limit Cycle Amplitude of the Motion Stage}}

\begin{figure}[!t]
    \centering
    \includegraphics[width = 0.6\linewidth]{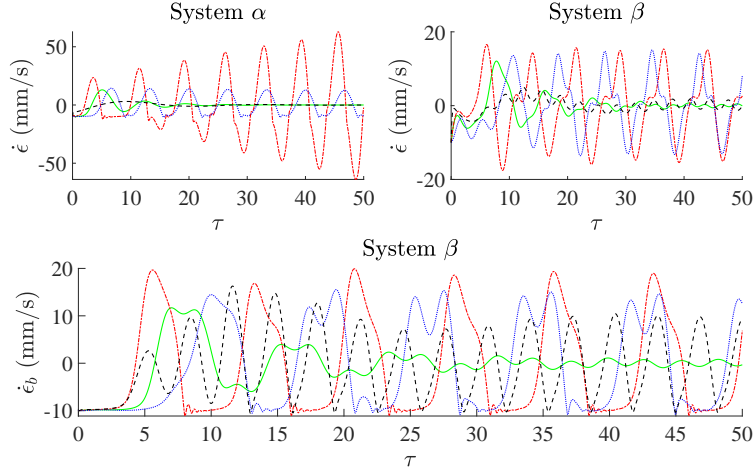}
    \caption{Comparison of velocity error time trajectories ($v_r=10$ mm/s) with different $k_i$ and $k_d$ for System $\alpha$ and $\beta$, where (1) green solid - default gains in Eq.(\ref{eq:standardGain}), (2) red dot - changes in $k_i=4.5 \times 10^6$ N/ms, (3) blue dash-dot - changes in $k_i=1 \times 10^5$ Ns/m and $k_d=20$ Ns/m, (4) black dash - changed to $k_d=1 \times 10^3$ Ns/m.}
    \label{fig:NumericalPID1}
\end{figure}

\begin{figure}[!b]
    \centering
    \includegraphics[width = 0.6\linewidth]{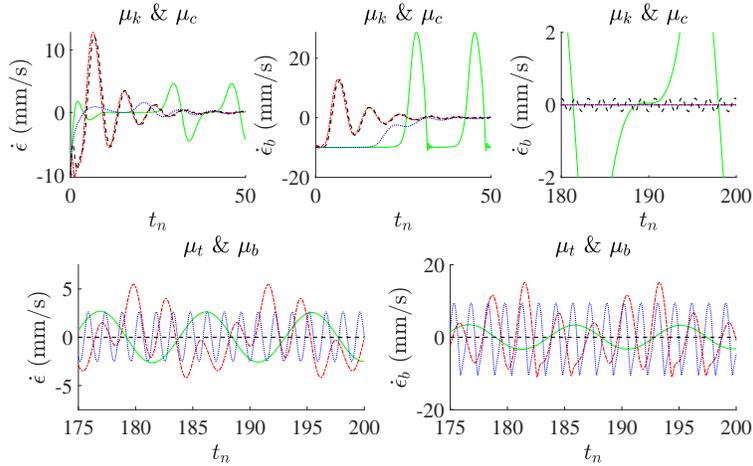}
    \caption{Comparison of velocity error time trajectories ($v_r = 10$ mm/s) with different scaling factors $\mu$ , where for $\mu_k \ \& \ \mu_c$: (1) green solid - $\mu_k = 0.1$, (2) red dot - $\mu_k = 10$, (3) blue dash-dot - $\mu_k = 0.1$ and $\mu_c = 100$, (4) black dash - $\mu_k = 3$ and $\mu_c = 0.25$; and for $\mu_t \ \& \ \mu_b$: (1) green solid - $\mu_b=2$, (2) red dot - $\mu_t = 2$, (3) blue dash-dot - $\mu_b = 0.5$, (4) black dash - $\mu_t = 0.5$.}
    \label{fig:NumericalParam1}
\end{figure}

\begin{figure}[!t]
    \centering
    \includegraphics[width = 0.6\linewidth]{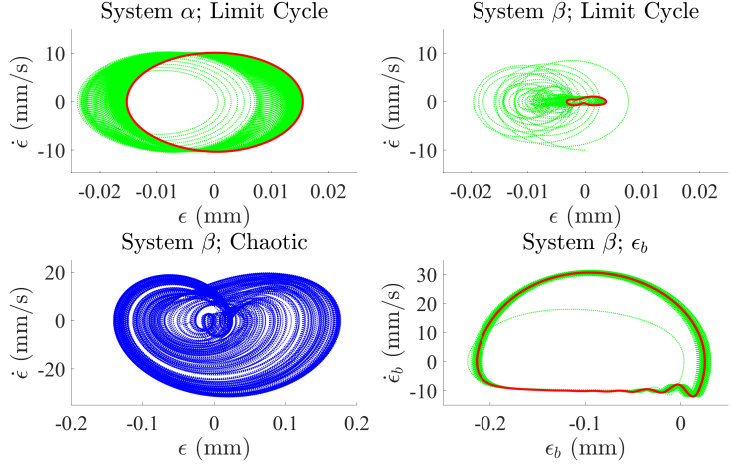}
    \caption{Phase portraits with limit cycles highlighted in solid red lines, where the control gains are (1) System $\alpha$, limit cycle  - $k_i= 1 \times 10^6$ N/ms, $k_p= 6.5 \times 10^5$ N/m, $k_d= 20$ Ns/m; (2) System $\beta$, limit cycle  - $k_i= 1 \times 10^6$ N/ms, $k_p= 6.5 \times 10^5$ N/m, $k_d= 20$ Ns/m; (3) System $\beta$, Chaotic  - $k_i= 8 \times 10^6$ N/ms, $k_p= 2 \times 10^4$ N/m, $k_d= 2 \times 10^2$ Ns/m; and (4) System $\beta$, $\epsilon_b$ - $k_i= 1 \times 10^7$ N/ms, $k_p= 2 \times 10^4$ N/m, $k_d= 2 \times 10^2$ Ns/m.}
    \label{fig: limitCycles}
\end{figure}

\begin{figure}[!b]
    \centering
    \includegraphics[width = 0.6\linewidth]{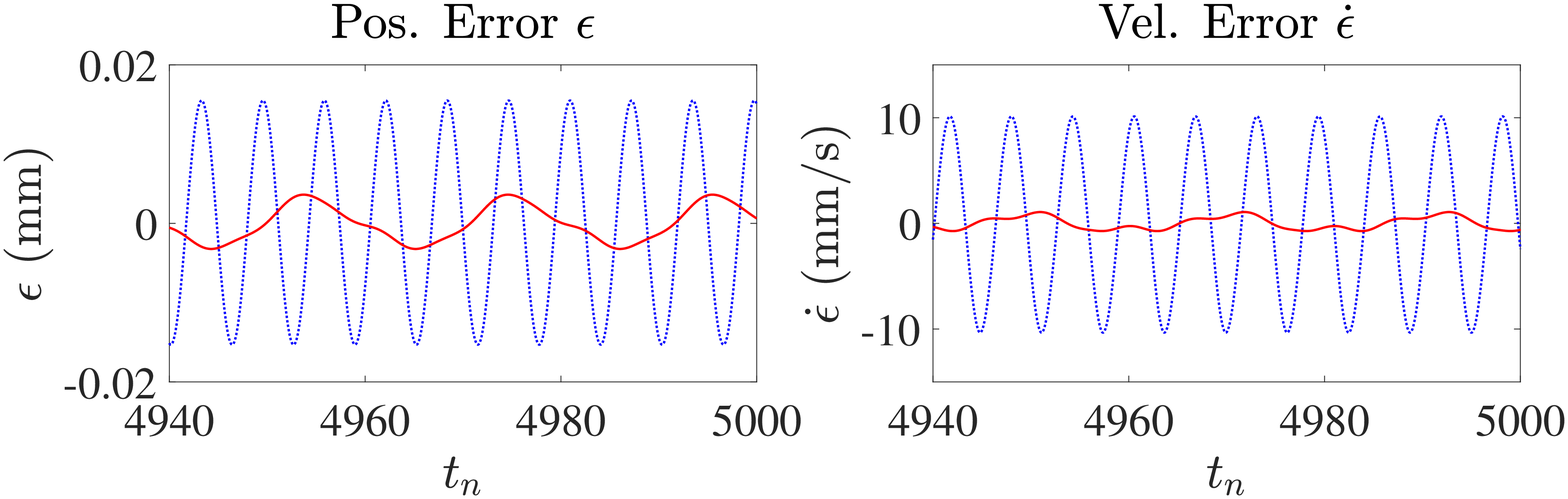}
    \caption{Comparison of limit cycle time trajectories, where - (1) blue dot: System $\alpha$, and (2) red solid: System $\beta$. The PID control gains for both systems are selected as $k_i= 1 \times 10^6$ N/ms, $k_p= 6.5 \times 10^5$ N/m, and $k_d= 20$ Ns/m.}
    \label{fig:amplitudeComparison}
\end{figure}

This subsection examines the role of the FI on the limit cycle amplitude of the motion stage using numerical simulations. Figure \ref{fig:NumericalPID1} shows the time domain data of $\Dot{\epsilon}$ and $\Dot{\epsilon}_b$ with different controller gains (i.e., $k_i$ and $k_d$). Note that non-dimensional time $t_n=\omega_n t$ is used to provide a good time span for observation. For both Systems $\alpha$ and $\beta$, cases (2) and (3) are unstable. However, two types of instability are observed for System $\alpha$: the instability in case (2) is caused by an extremely large $k_i$ that, eventually leads $\Dot{\epsilon}$ to infinity; the instability in case (3) is the stick-slip phenomenon caused by the initial instability of the PD controller as a result of the coupling between friction dynamics and multibody dynamics (i.e., refer to the stability region in Fig. \ref{fig:stabilityAreaPD2D}). Although the ranges of stable PID gains are reduced in the presence of the FI, case (2) indicates that the implementation of the FI may prevent the error from going unbounded. In addition, while a large $k_d$ causes instability for System $\beta$ as shown in case (4), it can be observed that the corresponding $\Dot{\epsilon}$ has very small oscillation amplitude which may not affect the precision of the motion stage in practice.

The simulated velocity errors using different design parameter scaling factors are shown in Fig. \ref{fig:NumericalParam1}. In general, the system stability with different combinations of design parameters is well predicted by the linear analysis. Interestingly, the unstable oscillation in case (4) does not show any pattern of stick-slip, even though the instability is related to the coupling between  friction dynamics and multibody dynamics. The numerical results from different $\mu_t$ and $\mu_b$ combinations in Fig. \ref{fig:NumericalParam1} match with the observation in Fig. \ref{fig:massStability}, which shows that the system stability is very sensitive to the change of masses in the presence of FI, since only case (4) is stable with the choice of scaling factors between 0.5 to 2. 

The above numerical results show that the majority of the instability behaviors are bounded, indicating the existence of limit cycles. Existing literature \pCite{johanastrom2008revisiting,saha2016modified} has studied the amplitudes of limit cycles in the friction oscillator with LuGre dynamics (which is equivalent to System $\alpha$). Therefore, phase portraits of the systems under different control gains are compared in Fig. \ref{fig: limitCycles}. In the limit cycle subfigures, the same parameters are used for Systems $\alpha$ and $\beta$. Note that  the limit cycle amplitude of System $\beta$ is significantly smaller than that of System $\alpha$.  The corresponding time domain data of $\epsilon$ and $\dot{\epsilon}$ are also plotted in Fig. \ref{fig:amplitudeComparison}. This shows that the introduction of FI can reduce the amplitude of friction-induced vibration, thus agreeing with the previous experimental studies \pCite{dong2017simple, dong2018experimental}. Sub-figure (3) illustrates that chaotic behavior may be observed with certain parameters in the presence of FI, and sub-figure (4) demonstrates the limit cycle of the bearing states $\epsilon_b$, where the stick-slip effect can be easily noticed.

The numerical study has presented many interesting observations about the characteristics of the systems. While these results are in agreement with the stability analysis, many nonlinear features of the system, such as the cause of chaotic behaviors and the reduction of limit cycle amplitudes with the FI, will need further study using nonlinear analysis.

%% file: sec_Conclusion.tex
\section{Conclusion and Future Work}
\label{sec:conclusion}

This paper analytically and numerically examined the influence of friction isolator on the dynamics of a PID controlled motion stage under the LuGre friction dynamics. Linear stability analysis was performed at the slipping equilibrium point of the systems.  The eigenvalues and stability of the system were parametrically studied with respect to the PID control gains, FI design parameters, and friction parameters. Then the numerical analysis was carried out, \bluenote{which validated the analytical results from linear stability analysis}, and provided further insights into the nonlinear behavior of the system. The main results are:
\begin{enumerate}[{(}1{)}]
    \item The effects of the friction parameters on the stability of the system with FI share similar characteristics as that without FI. Unless a very large Integral gain is used, the stability of the system under the PID controller is dominated by  the Proportional and Derivative gains. 
    \item FI can increase the stability region. Large $k_{\FI}$ or $c_{\FI}$ may lead to larger regions for stable PID gain selection, in particular, allowing larger stable $k_i$ to be paired with small $k_p$ and $k_d$ for faster steady state error convergence. Raising $c_{\FI}$ and lowering $k_i$ also allow a more flexible table-bearing mass ratio and larger payload capacity.
    \item The numerical examples show that FI can reduce the amplitudes of limit cycles and prevent unbounded error, hence improve the precision of the motion stage. 
\end{enumerate}

The findings in this work also lay the foundation for future investigations, which include the \bluenote{experimental validation of the dynamical analysis results}, the nonlinear analysis of the motion stage that features the nonlinear FI stiffness and the friction model, and the optimization of FI design parameters for better performance and stability of the motion stage.